# IMMEDIATE CALCULATION OF SOME POISSON TYPE INTEGRALS USING SUPERMATHEMATICS CIRCULAR EX-CENTRIC FUNCTIONS

Florentin SMARANDACHE[1]  &  Mircea  Eugen  SELARIU[2]
[1]Chair of Department of Math & Sciences, University of New Mexico-Gallup, USA
[2]Polytechnic University of Timişoara, Romania

## 0. ABSTRACT

This article presents two methods, in parallel, of solving more complex integrals, among which is the Poisson's integral, in order to emphasize the obvious advantages of a new method of integration, which uses the supermathematics circular ex-centric functions.

We will specially analyze the possibilities of easy passing/changing of the supermathematics circular ex-centric functions of a centric variable $\alpha$ to the same functions of ex-centric variable $\theta$. The angle $\alpha$ is the angle at the center point O(0,0), which represents the centric variable and $\theta$ is the angle at the ex-center E(k,ε), representing the ex-centric variable. These are the angles from which the points $W_1$ and $W_2$ are visible on the unity circle – resulted from the intersection of the unity/trigonometric circle with the revolving straight line d around the ex-centric E(k,ε) – from O and from E, respectively.

## KEYWORDS AND ABBREVIATIONS

  **C** - Centric, Circular,
  **CC** - Circular Centric,
  **E** - Ex-centric,
  **EC** - Ex-centric Circular,
  **F** - Function,
  **H** - Hyperbolic,
  **PI** - Poisson Integral,
  **M** - Mathematics,
  **CM** - Centric Mathematics,
  **EM** - Ex-centric Mathematics,
  **SM** – Supermathematics,
  **FSM** - **F** & **SM, FSM_EC**- **FSM** & **EC, FSM_EH**-**FSM** & **EH.**

## 1. INTRODUCTION

The discovery of the ex-centric mathematics (**EM**), as a vast extension of the centric common/mathematics (**CM**), which together form the SuperMathematics (**SM**), allows new simpler approaches, for resolving more complex integrals, among which we present (11) the Poisson integral (**PI**) [1]. To emphasize the new integration method, we will present, in parallel, the classic method of solving, only for **PI**, presented in [1] and the new method which utilizes **SM**'s ex-centric circular functions (**EC**) [2], [3], [4] .

The **SM-EC** functions, which will be in the center of our attention, are the radial ex-centric functions rex θ and Rex α and the derivatives ex-centric dex θ and Dex α, functions which are independent of the reference system selected.

The functions rex θ, of ex-centric variable θ, of the principal determination 1 and secondary 2, defined on the whole real axis for numeric ex-centricity $k^2 < 1$, and for $k^2 > 1$ exist only in the interval $\Im \in (\theta_i, \theta_f)$, in which $\theta_{f,i} = \pi + \varepsilon \pm \arcsin(1/k)$, $\alpha_{f,i} = \theta_{f,i} + \beta_{f,i}$ are

(1) $\quad rex_{1,2}\,\theta = rex_{1,2}(\theta, E(k,\varepsilon)) = -k.\cos(\theta - \varepsilon) \pm \sqrt{1 - k^2 \sin^2(\theta - \eta)}$,

where **E(k, ε)** is a pole, called **ex-center**, which divides the straight line d ($d = d^+ \cup d^-$), revolving around this point, in the positive semi straight line $d^+$, on which is situated the first principal determination $rex_1\,\theta$, as function of ex-centric variable θ and, respectively, Rex $\alpha_1$, of centric variable α of the function and in the negative semi straight line $d^-$, on which is situated, along it, the second determination, secondary, of the function $rex_2\,\theta$ and Rex $\alpha_2$. The expressions of the same entities (1), as functions of centric variable α, which exist on the whole real axis, no matter which is the numeric ex-centricity k, are

(2) $\quad\quad\quad\quad Rex\,\alpha_{1,2} = \pm\sqrt{1 + k^2 - 2k.\cos(\alpha_{1,2} - \varepsilon)}$

These functions represent, as Prof. Dr. Math. Octav Gheorghiu observed, the **distance** in plane, as oriented segments, in polar coordinates, between two points: the ex-center E(k, ε) and the intersection points $W_{1,2}$ (1, $\alpha_{1,2}$) – between the straight line d and the unity circle CT[1,O(0,0)] with the center in the origin O of the system of coordinates axis, right Cartesian or polar reference point.

For an E which is interior to the unity disc, the segment $EW_1$ is situated on the positive direction of the semi straight line $d^+$, being, in this case, positive, that is Rex $\alpha_1$ = $rex_1\theta > 0$, while the oriented segment $EW_2$, positive oriented on the negative semi straight line is negative, that is Rex $\alpha_2$ = $rex_2\,\theta < 0$, as it can be seen in the Figure 1.

For $k = \pm 1$ at $\alpha_1 \in (0, 2\pi) \Rightarrow \theta \in (0, \pi)$ and at $\alpha_2 \in (0, 2\pi) \Rightarrow \theta \in (\pi, 2\pi)$. In other words, if, the straight line d rotates around $E(k,\varepsilon) \subset C(1,O)$ with an angular speed Ω ($\theta = \Omega t$), the points $W_{1,2}$ rotate on the unity circle C(1,O) with a double angular speed ($\alpha_{1,2} = 2\Omega t$) in half of the period and is stationary in the second half of the period ($\alpha_{1,2} = 0$) taking turns in E(k, ε).

If E is exterior to the unity disk, that is $|k| > 1$, then both determinations will be on the same semi straight line, being, successive, both positive and then, after the rotation of d of π, both negative, therefore are of the same sign, and this will make their product, in this case, positive, and, while in the precedent case, the product of the two determinations of the function was always negative (see Fig.1).

We must observe also that at k > 1 and for $\alpha_{1,2} \in (0, 2\pi) \Rightarrow \theta \in (\theta_i, \theta_f)$; the ex-centric variable θ diminishes the interval of existence of FSM-CE, between an initial value $\theta_i$ and a final one $\theta_f$, with as much as the numeric growth of ex-centricity k. For $k \to \infty$ the interval is reduced to a single point on the real axis **R**, for each determination.

The results presented so far, will also be obtained from the relations that will follow.

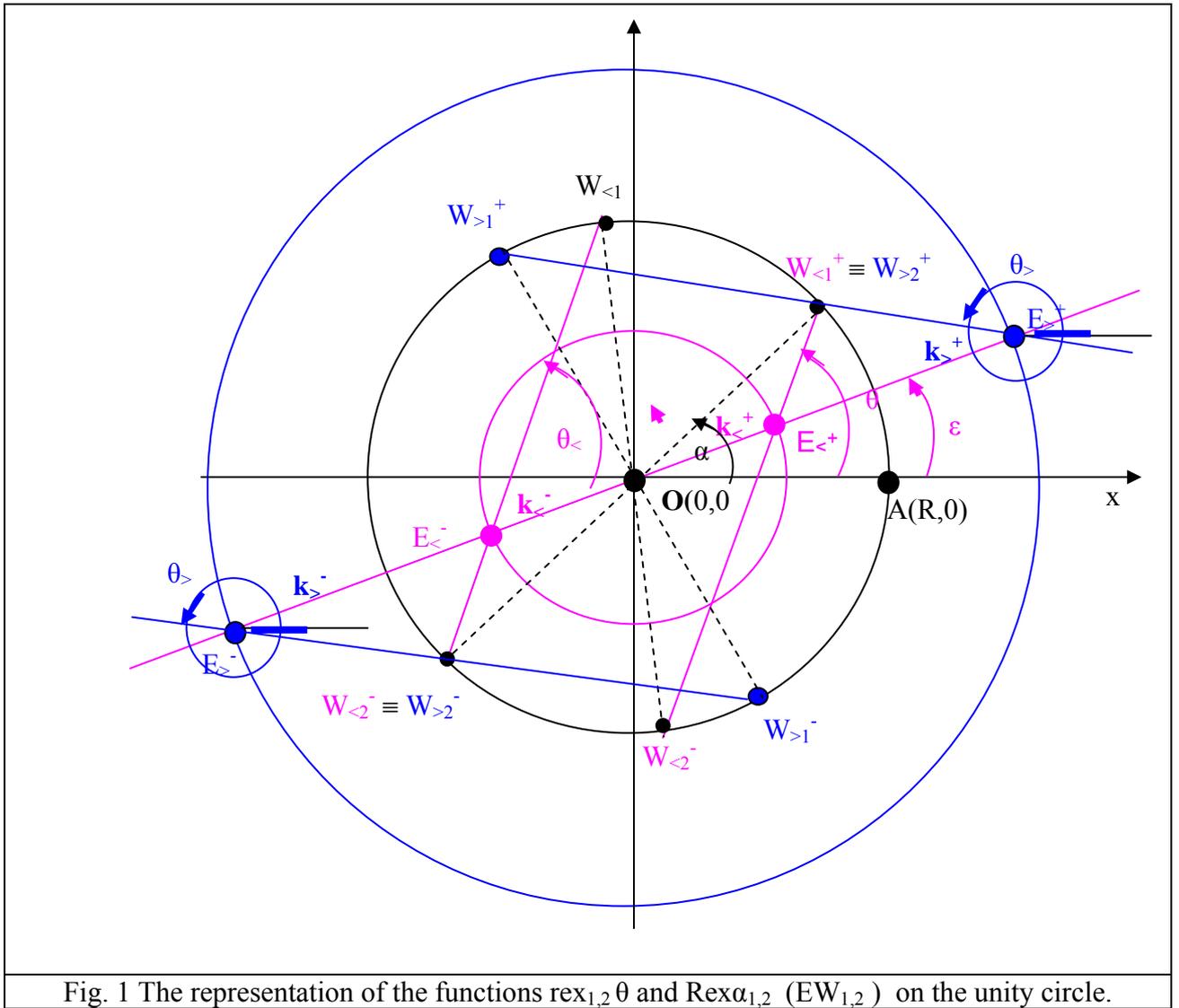

Fig. 1 The representation of the functions $rex_{1,2}\,\theta$ and $Rex\alpha_{1,2}$ ($EW_{1,2}$) on the unity circle.

The dependency between the two variables is:

(3) $\qquad \alpha_{1,2}(\theta) = \theta - \beta_{1,2}(\theta) = \theta \mp \arcsin[e.\sin(\theta - \varepsilon)]$

and, respectively

(4) $\quad \theta(\alpha_{1,2}) = \alpha_{1,2} + \beta(\alpha_{1,2}) = \alpha_{1,2} + \arcsin\left( \dfrac{k.\sin(\alpha_{1,2} - \varepsilon)}{\pm\sqrt{1 + k^2 - 2k.\cos(\alpha_{1,2} - \varepsilon)}} \right) =$

$\qquad = \alpha_{1,2} + \arcsin\left( \dfrac{k.\sin(\alpha_{1,2} - \varepsilon)}{\operatorname{Re} x\alpha_{1,2}} \right)$

or

(4') $\qquad \theta(\alpha_{1,2}) = \alpha_{1,2} + \arctan\left( \dfrac{k.\sin(\alpha_{1,2} - \varepsilon)}{1 - k.\cos(\alpha_{1,2} - \varepsilon)} \right) = \alpha_{1,2} + \arctan\left( \dfrac{\sin\beta(\alpha_{1,2})}{\cos\beta(\alpha_{1,2})} \right),$

where

(5) $\qquad \cos\beta(\alpha_{1,2}) = \dfrac{1 - k.\cos(\alpha_{1,2} - \varepsilon)}{\operatorname{Re} x\alpha_{1,2}}$

and

(6) $\quad \sin \beta (\alpha_{1,2}) = \dfrac{k.\sin(\alpha_{1,2} - \varepsilon)}{\operatorname{Rex}\alpha_{1,2}}$,

and the derivative of $d[\beta(\alpha)]/d\alpha$ is

(7) $\quad \dfrac{d\beta(\alpha)}{d\alpha} = \dfrac{k[\cos(\alpha - \varepsilon) - k]}{1 + k^2 - 2.k.\cos(\alpha - \varepsilon)} = \dfrac{k[\cos(\alpha - \varepsilon) - k]}{\operatorname{Rex}^2 \alpha_{1,2}}$

From (1), it results, without difficulty, that the sum, the difference, the product, and the ratio of the two determinations of the functions rex are:

(8) $\quad \begin{cases} \sum^{+} = rex_1\theta + rex_2\theta = -2k.\cos(\theta - \varepsilon) \\ \sum^{-} = rex_1\theta - rex_2\theta = 2\sqrt{1 - k^2.\sin^2(\theta - \varepsilon)} \\ \prod = rex_1\theta.rex_2\theta = \begin{cases} k^2 - 1 < 0, \rightarrow k < 1 \\ 0, \Rightarrow k = 1 \\ \pm(1 - k^2), \rightarrow k > 1 \end{cases} \\ \bigcup = \dfrac{|rex_2\theta|}{rex_1\theta} = \dfrac{d\alpha_2}{d\alpha_1} = \dfrac{|1 - k^2|}{\operatorname{Rex}^2\alpha_1} = \dfrac{d(rex_2\theta)}{d(rex_1\theta)} \end{cases}$

A function similarly useful, in this article, is the ex-centric derivative function of a centric variable $\alpha$, for which the form of expression is invariable at the position of the ex-center E is:

(9) $\quad \text{Dex } \alpha_{1,2} = \dfrac{1 - k.\cos(\alpha_{1,2} - \varepsilon)}{1 + k^2 - 2k.\cos(\alpha_{1,2} - \varepsilon)} = : \dfrac{d\theta}{d\alpha_{1,2}} = \dfrac{d(\alpha_{1,2} + \beta_{1,2})}{d\alpha} = \dfrac{1}{dex_{1,2}\theta}$,

and

(10) $\quad dex_{1,2} = 1 - \dfrac{k.\cos(\theta - \varepsilon)}{\sqrt{1 - k^2 \sin^2(\theta - \varepsilon)}}$

and the nucleus of Poisson integral

(11) $\quad \text{Nip } \alpha_{1,2} = \dfrac{d\gamma}{d\alpha_{1,2}} = \dfrac{d(\theta + \beta)}{d\alpha_{1,2}} = 1 + 2\dfrac{d\beta}{d\alpha_{1,2}} = \dfrac{1 - k^2}{1 + k^2 - 2k(\alpha_{1,2} - \varepsilon)} = \dfrac{1 - k^2}{\operatorname{Rex}^2\alpha_{1,2}} =$

$\quad = -\dfrac{\operatorname{Rex}\alpha_2}{\operatorname{Rex}\alpha_1}$

## 2. THE INTEGRATION USING THE CLASSIC METHOD [1]

The Poisson's integral, with modified notations, in accordance with the supermathematics ex-centric circular functions (SM - EC), is

(12) $\quad \text{PI } (k, \varepsilon) = \displaystyle\int_{-\pi}^{\pi} \dfrac{d\alpha}{1 + k^2 - 2.k.\cos(\alpha - \varepsilon)}$,

in which $k \in \Re$ and $\varepsilon \in [-\pi, \pi]$ are the parameters and, in the same time, the polar coordinates of the ex-center E. This is resolved in [1] as a simple integral which is dependent of a real parameter $\lambda \equiv k$, which will be further denoted as k, and representing, in EM, numeric ex-centricity k = e/R, the ratio between the real ex-centricity e and the circle radius R on which are placed the intersection points $W_1$ and $W_2$. The integral is simple, but the integration is quite laborious, as we will see later, and it will become indeed simple, only when passed from **CM** to **EM** with the utilization of the new supermathematics functions.

**Classical Solution**: The periodic real function

(13) $\quad f(\alpha) = \dfrac{1}{1+k^2 - 2k\cos(\alpha - \varepsilon)}$

is, as it is easily observed, the square of the radial ex-centric function of $\alpha$

(14) $\quad f(\alpha) = 1 / (\text{Rex}^2 \alpha)$,

defined for any $k \in \Re - \{\pm 1\}$ and $\varphi \in [-\pi, \pi]$.
**Remark:** Only one from the two determinations of the function Rex $\alpha_{1,2}$ is null (!) when E belongs to unity circle, that is $/k/ = 1$; the second determination having the expression which will be presented bellow. Based on the new knowledge from **EM**, now we can assert that the radial ex-centric function is defined also for $k = \pm 1$.
If $k = +1$, then

(15) $\quad \text{rex}_{1,2}\,\theta = -\cos(\theta-\varepsilon) \pm \sqrt{1 - \sin^2(\theta - \varepsilon)} \to \text{rex}_1\,\theta = \text{Rex}\,\alpha_1 = 0$

and $\text{rex}_2\theta = \text{Rex}\,\alpha_2 = -2\cos(\theta - \varepsilon)$ and, for $k = -1$, it results

(16) $\quad \text{rex}_{1,2}\,\theta = \cos(\theta-\varepsilon) \pm \sqrt{1 - \sin^2(\theta - \varepsilon)}$

such that, now, $\text{rex}_1\theta = \text{Rex}\,\alpha_1 = 2\cos(\theta-\varepsilon)$ and $\text{rex}_2\theta = \text{Rex}\,\alpha_2 = 0$, which it results and it can be seen, equally easy, also from the graphic.
Because

(17) $\quad \text{Rex}^2\,\alpha = [k - e^{i(\alpha - \varepsilon)}] \cdot [k - e^{-i(\alpha - \varepsilon)}] = [k - \text{rad}(\alpha-\varepsilon)] \cdot [k - \text{rad}-(\alpha-\varepsilon)]$,

in which, the radial centric functions [5], or in short, radial (denoted **rad**), equivalent to the exponential functions are unitary vectors, of symmetric directions, in relation to the straight line which contains the points O and E, therefore:

(18) $\quad \text{rad}(\alpha-\varepsilon) - \text{rad}[-(\alpha-\varepsilon)] = 2\cos(\alpha - \varepsilon)$

and

(19) $\quad \text{rad}(\alpha-\varepsilon) \cdot \text{rad}-(\alpha-\varepsilon) = \dfrac{\text{rad}(\alpha - \varepsilon)}{\text{rad}(\alpha - \varepsilon)} = 1$,

in which

(20) $\quad \text{rad}\,\alpha = e^{i\alpha}$

is equivalent, <u>in centric</u> (for $k = 0$, when $\alpha_1 = \theta$ and $\alpha_2 = \theta + \pi$) of functions rex $\theta$ and Rex$\alpha$ [5].
The function Rex$^2\alpha$ (16) has the roots:

(21) $\quad e^{\pm(\alpha - \varepsilon)} = \text{rad}[\pm(\alpha-\varepsilon)]$

which, for $\alpha = \varepsilon$ and also for $\alpha = \varepsilon - \pi$, din (14) it results

(22) $\quad k = \pm 1$.

By introducing in **PI** the variable $\alpha' = \alpha + \pi$, the change will lead to the integral:

(23) $\quad\quad\mathbf{PI}(-k) = \int_0^{2\pi} \dfrac{d\alpha'}{1 + k^2 + 2k.\cos(\alpha' + \varepsilon)},$

in which the numeric ex-centricity changes the sign, that is $k \rightarrow -k$, which is equivalent to the rotation of the ex-center E(k, ε) around the origin O (0,0) with π, on the circle with the radius k, that is $\varepsilon \rightarrow \pm (\varepsilon \pm \pi)$, or, yet, because of the inter-conversion properties of α with ε in the function cosine from (12), $\alpha \rightarrow \pm (\alpha \pm \pi)$.

Suppose that $k \neq \pm 1$, the change of the variable $\alpha' = \alpha + \pi$

(24) $\quad\quad z = e^{i(\alpha' + \varepsilon)},$

for which

(25) $\quad\quad dz / z = \dfrac{der(\alpha' + \varepsilon)d\alpha'}{rad(\alpha' + \varepsilon)} = \dfrac{i.rad(\alpha' + \varepsilon)d\alpha'}{rad(\alpha' + \varepsilon)} = i\, d\alpha'$

it will transform the segment [−π, + π] in the unity circumference, going in trigonometric positive sense ( sinistrorum / levogin). Then:

(26) $\quad\quad \mathbf{PI}(k, \alpha) = i \int_C \dfrac{dz}{(1+k^2)z + (1+z^2)k} = \dfrac{i}{k} \int_C \dfrac{dz}{z^2 + mz + 1},$

in which m = k +1/k.

The poly-functions f(z) from under the sign of $\int_C$ are z' = − k and z" = − 1 / k with the residues $a'_{-1}$ = Rez[f(z), -k] = k / (k² - 1) and $a''_{-1}$ = Rez[f(z), -1 / k] = k / (1 − k²), such that $a'_{-1} + a''_{-1} = 0$.

By applying the residues and semi-residues theorems, it results that for any angle ε ∈ [ −π, +π ]

(27) $\quad\quad \mathbf{IP}(k, \varepsilon) = \begin{cases} \dfrac{2\pi}{1-k^2}, & \text{for } |k| < 1; \\ 0, & \text{for } k = \pm 1 \\ \dfrac{2\pi}{k^2 - 1}, & \text{for } |k| > 1 \end{cases}$

The zero value for $\lim\limits_{k \rightarrow \pm 1} PI(k,\varepsilon)$ can be found choosing the contour Γ made of the circumferences C and γ (Fig. 2) , the last having the center in z" = 1 / k and the radius r < 1, from which we suppressed the interior portions of the reunion of the two circles. In this conditions, the integral $\int_\Gamma \dfrac{dz}{z^2 - kz + 1}$ is null even when k → 1 (or −1), appearing as a principal value in the Cauchy sense. We can then write:

(28) $$\mathbf{PI}(k, \varepsilon) = \begin{cases} \dfrac{2\pi}{\left|1-k^2\right|}, & \text{for} \Rightarrow k = \Re - \{\pm 1\} \\ 0, & \text{for} \Rightarrow k = \pm 1 \end{cases}$$

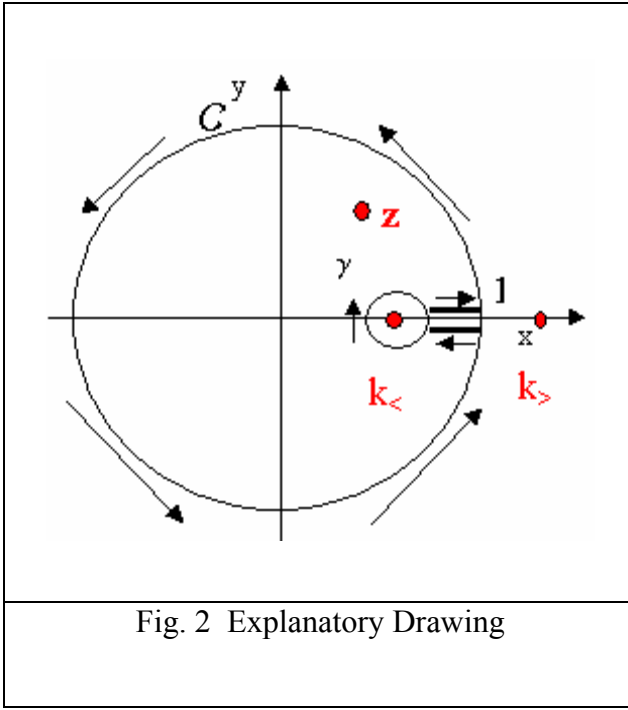
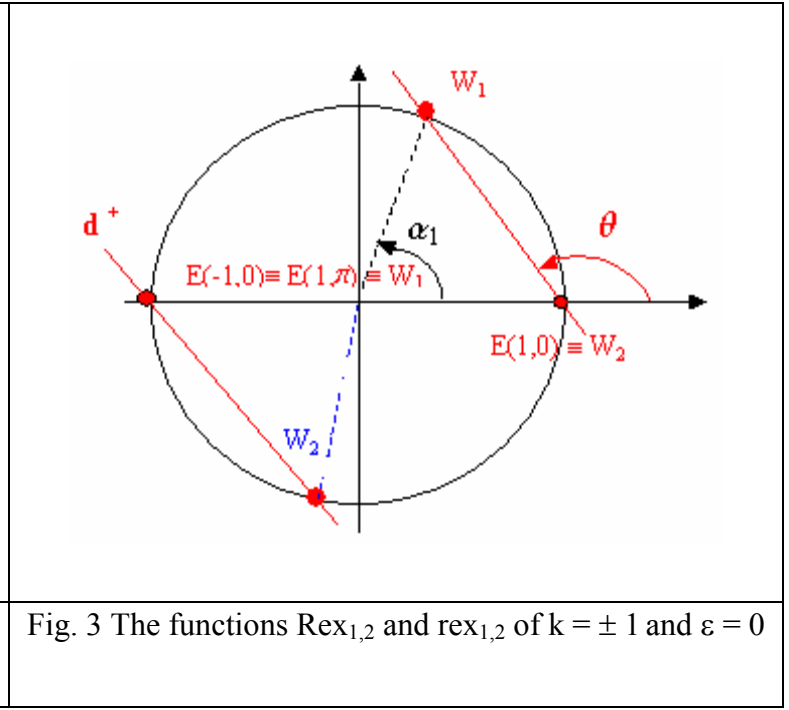

| Fig. 2 Explanatory Drawing | Fig. 3 The functions $Rex_{1,2}$ and $rex_{1,2}$ of $k = \pm 1$ and $\varepsilon = 0$ |

The result (28), presented in [1], can also be established directly, knowing that from (14), for $k = 1$, $rex_1\theta = 0$, and for $k = -1$, $rex_2\theta = 2\cos(\alpha-\varepsilon)$ such that:

(28) $$\int_{-\pi}^{\pi} \frac{d\alpha}{Rex^2\alpha_2} = \int_{-\pi}^{\pi} \frac{d\alpha}{4\cos^2(\alpha-\varepsilon)} = \frac{1}{4}\left|\tan(\alpha-\varepsilon)\right|_{-\pi}^{\pi} = \frac{1}{4}[\tan(\pi-\varepsilon) - \tan(-\pi-\varepsilon)] = 0$$

For $k \neq \pm 1$, will present the integral below.

### 3. THE INTEGRATION WITH THE HELP OF CIRCULAR EX-CENTRIC SUPERMATHEMATICS FUNCTIONS.

Multiplying $\mathbf{PI}(k, \varepsilon)$ with $(1-k^2) / (1-k^2)$ it results

(29) $$\mathbf{PI}(k, \varepsilon) = \frac{1}{1-k^2}\int_{-\pi}^{\pi} \frac{1-k^2}{Rex^2\alpha_{1,2}}d\alpha = \frac{1}{1-k^2}\int_{-\pi}^{\pi} -\frac{Rex_2\alpha}{Rex_1\alpha}d\alpha =$$
$$= \frac{1}{1-k^2}\int_{-\pi}^{\pi} \frac{d(\theta+\beta)}{d\alpha}d\alpha = \frac{1}{1-k^2}\int_{-\pi}^{\pi} d(\theta+\beta) = \frac{1}{1-k^2}\left|\theta+\beta\right|_{-\pi}^{\pi} = \frac{2\pi}{1-k^2},$$

for $k < 1$ and

(30) $$\mathbf{PI}(k, \varepsilon) = \frac{1}{1-k^2}\int_{-\pi}^{\pi} \frac{Rex\alpha_2}{Rex\alpha_1} = \frac{-2\pi}{1-k^2},$$

for k > 1, in which we took in consideration the relation (9) and the sign of the functions Rex $\alpha_{1,2}$, for k < 1 and for k > 1, that is an ex-center interior or exterior to unity disk and of the relation, for k < 1. The relations between the integration limits, taking in account of the dependencies [2]

(31) $\quad \begin{cases} \alpha_1 = \theta - \beta \\ \alpha_2 = \theta + \beta + \pi \end{cases}$

knowing that $\beta_1 + \beta_2 = \pi$ [3] are :

(32) $\quad$ If $\alpha_1 \in [-\pi, \pi] \Rightarrow$ then $\theta \in [\pi - \beta_1, \pi + \beta_1]$

and their difference is $+2\pi$, and if

(33) $\quad \alpha_2 \in [-\pi, \pi]$, then $\Rightarrow \theta \in [-2\pi - \beta, -\beta]$,

as it can be seen also in Fig. 1, and their difference is $-2\pi$.

(34) $\quad 1 + 2\dfrac{d\beta}{d\alpha} = 1 + 2\dfrac{e[(\cos(\alpha - \varepsilon) - e]}{\operatorname{Re} x^2 \alpha_{1,2}} = \dfrac{1 - k^2}{\operatorname{Re} x^2 \alpha_{1,2}} = \dfrac{d\alpha}{d\alpha} + \dfrac{d(2\beta)}{d\alpha} = \dfrac{d(\alpha + \beta + \beta)}{d\alpha} = \dfrac{d(\theta + \beta)}{d\alpha}$

because $\theta = \alpha + \beta$, and for

(35) $\quad \alpha = \begin{cases} \gamma_1 = -\pi \to \theta = -\pi + 2\beta_{1,2} \\ \gamma_2 = \pi \to \theta = \pi + 2\beta_{1,2} \end{cases}$,

as it results also from the figure, therefore

(36) $\quad \gamma_2 - \gamma_1 = 2\pi$

### CONCLUSIONS

Because of the labor volume in the two variants, the conclusion is, evidently, in the favor of the new method of integration, taking in account, firstly, the degree of complexity of the integration.

By utilizing the existing relations in EM, as, for example the relation (28), which can be written by denoting $\gamma = \theta + \beta$, from which $d\gamma = d(\theta + \beta)$, but $\alpha = (\theta - \beta)$ and $d\alpha = d(\gamma - 2\beta)$ or $d\alpha = d(\theta - \beta)$, such that

$$d\gamma / d\alpha = 1 + 2.d\beta / d\alpha = -\dfrac{k.\cos(\theta - \varepsilon)}{rex_1\theta} = -\dfrac{\operatorname{Re} x\alpha_2}{\operatorname{Re} x\alpha_1} = \dfrac{1 - k^2}{\operatorname{Re} x^2 \alpha_{1,2}}$$

and **PI** is an immediate integral, k being a constant parameter, as we saw before. Furthermore, from the relation (29) it results the Poisson's integral value undefined as being:

(37) $\quad \mathbf{IP_N} = \int \dfrac{d\alpha}{\operatorname{Re} x^2 \alpha} = \dfrac{1}{|1 - k^2|}[\theta(\alpha) + \beta(\alpha)] = \dfrac{1}{|1 - k^2|}[\alpha + 2\beta(\alpha)] = \dfrac{\alpha + 2\arcsin\dfrac{k\sin(\alpha\varepsilon)}{\operatorname{Re} x\alpha}}{|1 - k^2|}$

The integrals calculated in [1] with the help of the residues theorem

$$(38) \quad I_1 = \int_0^{2\pi} \frac{R - r\cos(\alpha - \varepsilon)}{R^2 + r^2 - 2r\cos(\alpha - \varepsilon)},$$

in which with r = k.R we denoted the real ex-centricity and with R the radius of a certain circle and

$$(39) \quad I_2 = \int_0^{2\pi} \frac{r\sin(\alpha - \varepsilon)}{R^2 + r^2 - 2r\cos(\alpha - \varepsilon)}$$

which, by the classical method presented in [1, pp. 186-187] are equally laborious and, unfortunately, wrong; by the new method, from EM, these integrals are immediate. Reducing R from (38) and (39) it results the functions to be integrated:

$$(40) \quad F_1 = \frac{R - r\cos(\alpha - \varepsilon)}{R^2 + r^2 - 2r\cos(\alpha - \varepsilon)} = \frac{1 - k\cos(\alpha - \varepsilon)}{\operatorname{Re} x^2 \alpha} = \text{Dex } \alpha_{1,2} = \frac{d\theta}{d\alpha}$$

such that the undefined integral is

$$(41) \quad I_{1N} = \int \frac{d\theta}{d\alpha} d\alpha = \int d\theta = \theta(\alpha_{1,2}) = \alpha + \arcsin\left[\frac{k\sin(\alpha_{1,2} - \varepsilon)}{\pm\sqrt{1 + k^2 - 2k\cos(\alpha_{1,2} - \varepsilon)}}\right],$$

such that, the defined integral (38) will be:

➢ For k = +1 ⇒ I$_1$ = π, because in the first determination 1 (principal) θ(α = 0) = π/2 and θ(α=2π) = 3π/2 and the difference is π. If k = − 1 for the first determination θ(α=0) = π and θ(α = 2π) = 2π, such that the difference is the same π. It results that for |k| = 1 ⇒ **I$_1$ = π.**
➢ For k > 1, the integral value I$_1$ is 0.

➢ For k < 1, the integral value is 2π.

The undefined integral I$_{2N}$ is:

$$(42) \quad I_{2N} = \int \frac{k.\sin(\alpha - \varepsilon)}{1 + k^2 - 2k\cos(\alpha - \varepsilon)} d\alpha = \int \frac{k.\sin(\alpha - \varepsilon)}{\operatorname{Re} x^2 \alpha} d\alpha = \int \frac{1}{\operatorname{Re} x\alpha} \frac{d(\operatorname{Re} x\alpha)}{d\alpha} d\alpha = \ln|\operatorname{Re} x\alpha|$$

Therefore, the defined integral I$_2$ is:

$$(43) \quad I_2 = \left|\ln\left|\operatorname{Re} x\alpha\right|\right|_0^{2\pi} = 0,$$

for any k and ε, knowing that Rex0 = rex0 = Rex2π = rex2π.

More integrals can be resolved immediately in this way without difficulties, if one knows the expressions of some supermathematics functions.
More integrals are presented in [6].